\RequirePackage{fix-cm}
\documentclass[smallextended, envcountsame]{svjour3}
\smartqed
\def\makeheadbox{\relax}
\journalname{}

\usepackage[unicode = false,
    pdftoolbar = true,
    colorlinks = true,
    linkcolor = blue,
    citecolor = blue,
    filecolor = black,
    urlcolor = blue,
    breaklinks = true]{hyperref}
\usepackage[ruled]{algorithm2e}
\usepackage{float}
\usepackage[utf8]{inputenc}
\usepackage{listings}           % For Code 
\usepackage{amsmath}            % For Math, duh
\usepackage{amssymb}            % For mathsfrs
\usepackage{xcolor}             % For Code Color
\usepackage[colorinlistoftodos]{todonotes}
\usepackage{graphicx}           
\usepackage[ruled]{algorithm2e}
\usepackage{verbatim}
\usepackage{geometry}
\usepackage{mcode}
\usepackage{caption}
\usepackage{enumerate}
\usepackage[thinlines, thiklines]{easybmat}
\usepackage{hhline}
\usepackage{nicematrix}
\usepackage{mathtools}
\usepackage[noblocks]{authblk}

\DeclarePairedDelimiter{\floor}{\lfloor}{\rfloor}
\DeclarePairedDelimiter{\ceil}{\lceil}{\rceil}
\newcommand{\R}{\mathbb{R}}
\newcommand{\Y}{\mathbb{Y}}
\newcommand{\D}{\mathbb{D}}
\newcommand{\B}{\mathbb{B}}
\newcommand{\G}{\mathbf{G}}
\newcommand{\A}{\mathbf{A}}
\newcommand{\V}{c\mathbf{V}}

\newcommand{\set}{\mathbb{S}}

\newcommand{\N}{\mathbb{N}}

\newcommand{\one}{\bold{1}}
\newcommand{\zero}{\bold{0}}

\newcommand{\p}{\mathring{p}_{\B_n}}

\newcommand{\bbm}{\begin{bmatrix}}
\newcommand{\ebm}{\end{bmatrix}}

\newcommand{\Pe}{\mathcal{P}}

\DeclareMathOperator{\rem}{rem}

\DeclareMathOperator{\pspan}{pspan}

\DeclareMathOperator{\von}{von}

\DeclareMathOperator{\cm}{cm}
\DeclareMathOperator{\rank}{rank}

\DeclareMathOperator{\spann}{span}

\DeclareMathOperator{\argminop}{argmin}

\DeclareMathOperator{\id}{Id}
\DeclareMathOperator{\diag}{Diag}

\DeclareMathOperator{\proj}{proj}

\DeclareMathOperator{\minimize}{Minimize}
\DeclareMathOperator{\st}{subject\, to}

\newcommand{\argmin}[1]{\underset{#1}{\argminop}}

\begin{document}

\title{\centering Nicely structured positive bases with maximal cosine measure}
\author{Warren Hare \and Gabriel Jarry-Bolduc \and Chayne Planiden}
\institute{Warren Hare \at
              Department of Mathematics, University of British Columbia, Kelowna, \\
              British Columbia, Canada.  Hare's research is partially funded by the Natural Sciences and Engineering Research Council (NSERC) of Canada, Discover Grant \#2018-03865. ORCID 0000-0002-4240-3903\\
              \email{warren.hare@ubc.ca}           %  \\
           \and
          Gabriel Jarry-Bolduc\at
           Department of Mathematics, University of British Columbia, Kelowna, \\
              British Columbia, Canada.  Jarry-Bolduc's research is partially funded through the Natural Sciences and Engineering Research Council (NSERC) of Canada, Discover Grant \#2018-03865. ORCID 0000-0002-1827-8508\\
              \email{gabjarry@alumni.ubc.ca} 
              \and
              Chayne Planiden \at
              School of Mathematics and Applied Statistics, University of Wollongong, Wollongong, NSW, 2500, Australia. Research supported by University of Wollongong. ORCID 0000-0002-0412-8445\\
              \email{chayne@uow.edu.au}
              }
\date{\today}

\maketitle
\begin{abstract}
The properties of positive bases make them a useful tool in derivative-free optimization (DFO) and an interesting concept in mathematics. The notion of the \emph{cosine measure} helps to quantify the quality of a positive basis. It provides information on how well the vectors in the positive basis uniformly cover the space considered. The  number of vectors in a positive basis is known to be between $n+1$ and $2n$ inclusively. When the number of vectors is strictly between $n+1$ and $2n$, we say that it is an intermediate positive basis. In this paper, the structure of intermediate positive bases with maximal cosine measure is investigated. The structure of an intermediate  positive basis with maximal cosine measure over a certain subset of positive bases is provided. This type of positive bases has a simple structure that makes them easy to generate with a computer software. 
\end{abstract}
\keywords{Positive basis \and Positive spanning set \and Cosine measure \and Derivative-free optimization}
\section{Introduction}\label{sec:intro}

The concept of positive bases was developed in 1954 by Davis \cite{Davis1954} and their properties were further developed in \cite{Mckinney1962,Romanowicz1987,Shephard1971}. A positive basis of a space is a set of vectors than spans the space using only positive coefficients in their linear combinations, as opposed to standard bases that allow negative coefficients as well. It is known that the number of vectors $s$  in a positive basis of $\R^n$ is between $n+1$ and $2n$ inclusive.  This result is proved in \cite{Davis1954} and an alternate short proof for the upper bound has been published in \cite{Audet2011}. A positive basis of size $n+1$ is called a minimal positive basis and one of size $2n$ is called a maximal positive basis. If $n+1<s<2n$, then the positive basis is known as a positive basis of intermediate size. The structure of minimal and maximal positive bases is now well-understood and has been meticulously characterized in \cite{Regis2016}. However, the structure of intermediate positive bases is not as obvious and few results are available in the literature. This paper explores intermediate positive bases.  It gives a definition of  optimal positive basis and provides a method for constructing an optimal positive basis over a set of positive bases with nice structure. One of the most important results on the structure of intermediate positive bases is provided in  \cite[Theorem 1]{Romanowicz1987}, which we use in proving the main results of this paper. 

In the last twenty years, the topic of positive bases and positive spanning sets has become an active area of research, due to its value in derivative-free optimization (DFO). As mentioned in \cite{Conn2009,Lewis1996}, the key property of positive bases is that they always  contain a descent direction for any function  whenever the gradient of the function exists and is non-zero at the point of interest. Positive bases are employed in direct search methods such as pattern search \cite{Audet2017,Custodio2008,Torczon1997,Vaz2009}, grid-based methods \cite{Coope2001,Coope2002} and many others \cite{Abramson2009,Audet2006,Audet2014,Kelley2011,Kolda2003,van2013using}. They have also be used in non-DFO methods.  For example, a feasible descent direction method that relies on the properties of positive spanning sets and minimizes a class of nonsmooth, nonconvex problems with linear constraints has been proposed in \cite{beck2020}.

The \emph{cosine measure} is used to indicate the quality of a positive basis \cite{Torczon1997} and it is this measure that determines whether the basis is optimal.  Roughly speaking, a high value of the cosine measure indicates that the positive basis covers the space more uniformly. In general, having higher cosine measure is preferable and convergence properties of certain DFO algorithms depend on the cosine measure of the positive basis employed. A deterministic algorithm to compute the cosine measure of any finite positive spanning set in finite time is introduced in \cite{hare2020cm}. Another method for computing the cosine measure of any finite set has also been introduced in \cite{regis2021}. 

In \cite{Naevdal2018}, the maximal cosine measures for maximal and minimal positive bases are found and the structure of positive bases attaining these upper bounds for the cosine measure are characterized.  However, the maximal cosine measures for intermediate positive bases was not developed.  In this paper, we investigate this issue in depth. We define two subsets of positive bases with nice properties.  We investigate the problem of finding an intermediate positive basis with maximal cosine measure on these two subsets.   We develop properties of this type of intermediate positive bases and show that the algorithm to compute the cosine measure introduced in \cite{hare2020cm} can be simplified in the presence of such positive bases.

This paper is organized as follows. In Section \ref{sec:prel}, the main definitions and background results necessary to understand this paper are presented. In Section \ref{sec:main}, we investigate the properties of intermediate positive bases with maximal cosine measure in $\R^3$. Then we consider a general space $\R^n$ in the second part of Section \ref{sec:main}. Finally, in Section \ref{sec:conc} we summarize the main results of this paper and discuss future research directions.

\section{Preliminaries} \label{sec:prel}
To begin, we present the definition and three important properties of a positive spanning set. The space considered throughout this paper is $\R^n.$
\begin{definition}[Positive span and positive spanning set of $\R^n$]
The \emph{positive span} of a set of vectors $\set=\begin{bmatrix}d_1& d_2&\cdots&  d_s \end{bmatrix}$ in $\R^n$, denoted $\pspan(\set)$, is the set $$\{v \in \R^n: v=\alpha_1d_1+ \dots+\alpha_s d_s, \alpha_i\geq 0, i=1, 2, \dots, s\}.$$
A \emph{positive spanning set} of $\R^n$  is a set of vectors $\set$ such that $\pspan(\set)=\R^n.$ 
\end{definition}

It will be convenient to regard a set of vectors as a matrix whose columns are the vectors in the set. The following three lemmas will be used in Section \ref{sec:main}.
%=========================================================================================================
\begin{lemma}\textup{\cite[Theorem 2.3]{Conn2009}} \label{lem:posspanlem}
Let $\set=\bbm d_1&\cdots&d_s \ebm$ be a set of vectors that spans $\R^n.$ Then $\set$ positively spans $\R^n$ if and only if there exist real scalars $\alpha_1, \dots, \alpha_s>0$ such that $\sum_{i=1}^s \alpha_id_i=0.$
\end{lemma}
 
%=====================================================================================================

\begin{lemma}\textup{\cite[Theorem 2.3]{Regis2016}} \label{lem:Removing}
If $\set=\begin{bmatrix}d_1&&\cdots&d_s \end{bmatrix}$ positively spans $\R^n$, then $\set \setminus\{d_i\}$ linearly spans $\R^n$ for any $i\in\{1, \dots, s\}$.
\end{lemma}
\begin{lemma}\cite[Theorem 2.3]{Conn2009} \label{lem:dotprodneg}
Let $\set=\bbm d_1&\cdots&d_s \ebm$ be a set of non-zero vectors in $\R^n$. Then the 
set $\set$ is a positive spanning set of $\R^n$ if and only if
for every non-zero vector $v$ in $\R^n$, there exists an index $i \in \{1,2, \dots, s\}$ such that $v^\top d_i<0.$
\end{lemma}

To define a positive basis of $\R^n$, we  introduce the concept of positive independence.

%Def: Pos ind.
\begin{definition}[Positive independence]
A set of vectors $\set=\begin{bmatrix}d_1& d_2&\cdots&  d_s \end{bmatrix}$ in $\R^n$  is \emph{positively independent}  if  $d_i \notin \pspan(\set\setminus {d_i})$ for all $i \in \{1, 2, \dots s\}.$
\end{definition}
%=====================================================================================================
%====================================================================================================
\begin{definition}[Positive basis of $\R^n$]\label{def:pbasis}
A \emph{positive basis} of $\R^n$ of size $s$, denoted $\D_{n,s}$, is a positively independent set whose positive span is $\R^n$.
\end{definition}
%=====================================================================================================

Equivalently, a positive basis of $\R^n$ can be defined as a set of non-zero vectors of $\R^n$ whose positive span is $\R^n,$ but for which no proper subset exhibits the same property \cite{Conn2009}.
For convenience, we will assume that the vectors in a  positive basis are unit vectors in this paper.

It is known that the size $s$ of a positive basis is  between $n+1$ and $2n$ inclusively.
We say that the positive basis $\D_{n,s}$ is \emph{intermediate size} if $n+1 <  s < 2n.$ Note that there is no positive basis of intermediate size when $n \in \{1,2\}.$ When $s=n+1,$ we say that the positive basis is minimal size. A positive basis of minimal size in $\R^n$ will be written $\D_n$ (the second argument in the subscript referring to the size may be omitted when $s=n+1$).  When $s=2n,$ we say that the positive basis is maximal size. If $s>n+1,$ we say the positive basis $\D_{n,s}$ is non-minimal.

The principal tool for determining the quality of a positive  basis and how well it covers the space $\R^n$ is the cosine measure.

\begin{definition}[Cosine measure]\label{def:cosinemeasure} 
The \emph{cosine measure} of a finite set $\set$ of non-zero vectors is defined by $$\cm(\set)=\min_{\substack{\Vert u \Vert=1\\u \, \in \, \R^n}} \max_{d \, \in \, \set} \frac{u^\top d}{\Vert d \Vert}.$$
\end{definition}
Values of the cosine measure near zero suggest the positive spanning property is approaching a deterioration. A high value of cosine measure indicates that the vectors in the set more uniformly cover the space; in other words, the vectors are spaced farther away from one another. 
Note that given any positive basis $\D_{n,s}$ where $n\geq 2$, the cosine measure is bounded by $0< \cm(\D_{n,s})<1$ \cite[Proposition 10]{hare2020cm}. 
When $n=1,$ there is only one positive basis of unit vectors:
\begin{align*}
    \D_{1,2}&=\bbm 1&-1\ebm.
\end{align*}
Its cosine measure is equal to 1. The positive basis $\D_{1,2}$ is both minimal size and maximal size.
\begin{definition}[Optimal positive basis]
 A positive basis $\D_{n,s}$ is \emph{optimal} over a non-empty set $\mathbb{P}$ if $$\cm(\D_{n,s})\geq \cm(\D'_{n,s})$$ for any positive basis $\D'_{n,s}$ in $\mathbb{P}.$
\end{definition}

The set of positive bases containing all positive bases of size $s$ in $\R^n$ will be denoted $\Pe.$ When we write that $\D_{n,s}$ is optimal without mentioning the set considered, it is implied that the set considered is  $\Pe.$ The structure and properties of optimal bases of minimal size  and maximal size  are well-known and have been rigorously proved in \cite{Naevdal2018}.  We will denote by $\mathring{\D}_n$ an optimal positive basis of minimal size. Note that in $\R,$ the only  positive basis of unit vectors $\D_{1,2}=\bbm 1&-1 \ebm$ is an optimal positive basis of minimal size and hence, will be denoted by $\mathring{\D}_1.$ 

Next, we define the cosine vector set, the active set of vectors and present a very important property of the active set.
%================================================================================================
\begin{definition}[The cosine vector set]
Let $\set$ be a set of non-zero vectors in $\R^n$. The cosine vector set of $\set$, denoted $\V(\set)$, is defined as $$\V(\set)=\argmin{\substack{\Vert u \Vert=1\\u \, \in \, \R^n}}\max_{d\, \in\, \set} \frac{u^\top d}{\Vert d \Vert}.$$ 
\end{definition}
%================================================================================================
%================================================================================================

\begin{definition}[The all activity set of vectors]
Let $\set$ be a set of unit vectors in $\R^n$ and let $u \in \V(\set).$ The \emph{active set} of vectors in $\set$ on $u$ denoted $\mathcal{A}(u,\set)$, is defined as $$\mathcal{A}(u,\set)=\left \{ d  \in \set:d^\top u=\cm(\set) \right \}.$$
The \emph{all activity set} is defined as $$\A(\set)=\bigcup_{u \in \V(\set)} \mathcal{A}(u,\set).$$
\end{definition}

\begin{proposition} \label{prop:activesetproperty}
Let $\D_{n,s}$ be a positive basis of $\R^n.$ Then the all activity set $\A(\D_{n,s})$ contains a basis of $\R^n.$
\end{proposition}
\begin{proof}
This follows directly from \cite[Corollary 18]{hare2020cm} as $\mathcal{A}(u,\D_{n,s})$ contains a basis of $\R^n$ for each $u$ in $\V(\D_{n,s}).$ \qed
\end{proof}

The next step is to introduce a deterministic algorithm to compute the cosine measure of any positive basis. First, let us recall the definition of a Gram matrix.
\begin{definition}[Gram matrix]
Let $\set=\begin{bmatrix}d_1&d_2& \cdots& d_s\end{bmatrix}$ be a set of vectors in $\R^n$ with dot product $d_i^\top d_j$. The Gram matrix of the vectors $d_1, d_2, \dots, d_s$ with respect to the dot product, denoted $\G(\set),$ is given by
$\G(\set)=\set^\top \set.$
\end{definition}

%\begin{lemma}\cite[Lemma 1]{Naevdal2018} \label{Lem:gammab} Let $\B=\begin{bmatrix}d_1&d_2&\cdots& d_n\end{bmatrix}$ be a basis of unit vectors in $\R^n$. Let $\one \in \R^n$ be the vector having all its entries equal to one. Then there exists a unit vector $u_\B \in \R^n$ such that $u_\B^\top d_i=\gamma_\B>0$ for all $i \in \{1, 2, \dots, n\}$ where $$\gamma_\B=\frac{1}{\sqrt{\one^\top \G(\B)^{-1}\one}}.$$
%Furthermore, $\gamma\leq \sqrt{e^TG(v_1, \dots, v_n)e}/n$.
%\end{lemma}
%===========================================================================================================

%Note that the unit vector $u_\B$ such that $u_\B^\top d_i=\gamma_{\B}$ for all $i$'s is unique since $\{d_1, \dots, d_n\}$ is a set of $n$ linearly independent vectors. Also, note that $\gamma_\B<1$ whenever $n\geq 2.$
%In fact, the positive value of the $n$ equal dot products is unique whenever $u$ is a unit vector and $\{d_1, \dots, d_n\}$ is a set of linearly independent vectors.

%============================================================================================================
%\begin{lemma}\cite[Lemma 13]{hare2020cm} \label{lem:gammab}
%Let $\B=\begin{bmatrix}d_1&d_2&\cdots& d_n\end{bmatrix}$ be a basis of unit vectors in $\R^n$. Suppose $u$ is a unit vector such that $u^\top d_1=\dots=u^\top d_n=\alpha>0.$ Then $\alpha=\gamma_\B,$ where $\gamma_\B$ is defined in Lemma \ref{Lem:gammab}.
%\end{lemma}

Algorithm \ref{algo:computecm} below was introduced in \cite[Algorithm 1]{hare2020cm}. 
\begin{center}
\scalebox{1}{
\begin{algorithm}[H] \label{algo:computecm}
\DontPrintSemicolon
\caption{The cosine measure of a positive basis  in $\R^n$ \label{alg:cmpbasis}}
Given $\D_{n,s}$, a positive basis of size $s$ in $\R^n$:

         \textbf{1. For all bases $\B_n\subset \D_{n,s}$, compute\;  }
         {
  $ \displaystyle
    \begin{aligned}
      (1.1)&\quad\gamma_{\B_n}=\frac{1}{\sqrt{\one^\top \G(\B_n)^{-1}\one}} &&\text{(The positive value of the $n$ equal dot products)},\\
      (1.2)&\quad u_{\B_n}=\gamma_{\B_n}\B_n^{-\top} \one &&\text{(The unit vector associated to $\gamma_{\B_n}$)},\\
      (1.3)&\quad p_{\B_n}=\begin{bmatrix}p_{\B_n}^1&\cdots&p_{\B_n}^s \end{bmatrix}=u_{\B_n}^\top\D_{n,s} &&\text{(The dot product vector)}, \\
      (1.4)&\quad\p =\max_{1\leq i \leq s }p_{\B_n}^i &&\text{(The maximum value in $p_{\B_n}$)}.
    \end{aligned}
  $
\par}

        \textbf{2. Return  \;}
    	
    	 {
  $ \displaystyle
    \begin{aligned}
        (2.1)&\quad \cm(\D_{n,s})=\min_{\B_n \, \subset \, \D_{n,s}} \p &&\text{(The cosine measure of $\D_{n,s}$)}\\
        (2.2)&\quad\V(\D_{n,s})=\{u_{\B_n}:\p=\cm(\D_{n,s})\}&&\text{(The cosine vector set of $\D_{n,s}$)}.
    \end{aligned}
  $
\par}
\end{algorithm}}
\end{center}
In \cite{hare2020cm}, it was proved that if $\B_n$ is a basis of $\R^n$ contained in $\mathcal{A}(u,\D_{n,s})$ for some $u \in \V(\D_{n,s}),$ then the cosine measure of the positive basis $\D_{n,s}$ is given by
\begin{align}\label{eq:computecmwithbn}
    \cm(\D_{n,s})&=\frac{1}{\sqrt{\one^\top \G(\B_n)^{-1} \one }}.
\end{align}
Therefore, an optimal positive basis has a structure that maximizes \eqref{eq:computecmwithbn}, or equivalently, that minimizes $\one^\top \G(\B_n)^{-1} \one$.
The next lemma shows that  any basis of $\R^n$ contained in an optimal \textbf{minimal} positive basis $\mathring{\D}_n$ is in the all activity set.
\begin{lemma} \label{lem:optimalminimalDns}
Let $\mathring{\D}_{n}$ be an optimal minimal positive basis in $\R^n.$ Let $\B_n$ be any  basis of $\R^n$ contained in $\mathring{\D}_n.$ Then $\B_n$ is in $\A(\mathring{\D}_n)$
\end{lemma}
\begin{proof}
This follows from \cite[Theorem 1]{Naevdal2018} as all Gram matrices of a basis $\B_n$ contained in  $\mathring{\D}_n$  are equal. It follows from Lemma \ref{lem:dotprodneg} that $\gamma_{\B_n}=\mathring{p}_{\B_n}$ for all bases $\B_n$ (where $\mathring{p}_{\B_n}$ is defined as in Algorithm \ref{algo:computecm}). Therefore, any basis $\B_n$ contained in $\mathring{\D}_n$ is in $\A(\mathring{\D}_n).$  \qed
\end{proof}

Now we recall the definition of a principal submatrix and properties of positive definite matrices and positive semidefinite matrices.
\begin{definition}[Principal submatrix]\textup{\cite[ Section 0.7.1]{Horn1990}}
Let $A \in \R^{n \times m}.$ For index sets $\alpha \subseteq \{1, \dots, n\}$ and $\beta \subseteq \{1, \dots, m\},$ we denote by $A[\alpha, \beta]$ the submatrix of entries that lie in the rows of $A$ indexed by $\alpha$ and the columns indexed by $\beta.$ If $\alpha=\beta,$ the submatrix $A[\alpha,\alpha]$ is a principal submatrix of $A.$
\end{definition}
\begin{lemma}\textup{\cite[Obs. 7.1.2]{Horn1990}} \label{lem:principalsubmatrix}
Let $A \in \R^{n \times n}$ be a real symmetric matrix. If $A$ is positive definite, then all of its principal submatrices are positive definite.
\end{lemma}
\begin{lemma}\textup{\cite[Theorem 7.2.7]{Horn1990}} \label{lem:pdisgram}
Let $A$ be a symmetric matrix in $\R^{n \times n}$. The matrix $A$ is positive definite if and only if there is a $B \in \R^{m \times n}$ with full column rank such that $A=B^\top B.$ 
\end{lemma}
%\begin{lemma}\label{lem:uniquesquareroot}
%A positive definite matrix $A \in \R^{n \times n}$ has precisely one positive definite matrix $B \in \R^{n \times n}$ such that $A=B^\top B.$ \textcolor{blue}{Find a reference in a book. From Wikipedia.}
%\end{lemma}
\begin{lemma}\textup{\cite[Theorem 7.2.1]{Horn1990}}\label{lem:inversepdispd}
A nonsingular symmetric matrix $A \in \R^{n \times n}$ is positive definite if and only if $A^{-1}$ is positive definite.
\end{lemma}
%\begin{lemma}\textup{\cite[Corollary 7.7.4]{Horn1990}} \label{lem:pdinverseineq}
%Let $A,B \in \R^{n \times n}$ be symmetric matrices. If $A$  and $B$ are positive definite matrices, then $A-B$ is positive semi-definite if and only if $B^{-1}-A^{-1}$ is positive semi-definite.
%\end{lemma}

Note that Lemma \ref{lem:pdisgram} and Lemma \ref{lem:inversepdispd} explain  why the value $\gamma_{\B_n}$ defined in Algorithm \ref{algo:computecm} Step (1.1)  is  positive for any basis $\B_n$ of $\R^n.$
%\begin{lemma}\textup{\cite[Observation 7.1.8]{Horn1990}} \label{lem:CAC}
%Let $A \in \R^{n \times n}$ be symmetric and let $C \in \R^{n \times m}.$ 
%\begin{enumerate}[(a)]
%\item If $A$ is positive semidefinite, then $C^\top A C$ is positive semidefinite.
%\item If $A$ is positive definite and $\rank(C)=m,$ then $C^\top A C$ is positive definite.
%\end{enumerate}
%\end{lemma}
The last lemma of this section will  be helpful to find the cosine measure of a positive basis $\D_{n,s}$ when  a basis $\B_n \in \A(u,\D_{n,s})$ is written as a block matrix. 
\begin{lemma}[Inverse of a block matrix]\textup{\cite[Corollary 4.1]{Lu2002}} \label{lem:inverseblock}
Let $A$ be a symmetric invertible matrix and $D$ be a symmetric matrix.
The inverse of a positive definite matrix having the form $$G=\begin{bmatrix}A&B^\top\\B&D\end{bmatrix}$$ is $$G^{-1}=\begin{bmatrix}A^{-1}+A^{-1}B^\top(D-BA^{-1}B^\top)^{-1}BA^{-1}&\quad -A^{-1}B^\top(D-BA^{-1}B^\top)^{-1}\\-(D-BA^{-1}B^\top)^{-1}BA^{-1}&(D-BA^{-1}B^\top)^{-1} \end{bmatrix}.$$
\end{lemma}

%=====================================================
%Deleted Section 

\section{Main Results} \label{sec:main}
In 1987, Romanowicz introduced the concepts of sub-basis, sub-positive basis and critical vectors \cite{Romanowicz1987}. He used these concepts to characterize the structure of positive bases. As such, they will be helpful to characterize the structure of any  non-minimal positive basis.

\begin{definition}[Sub-basis and sub-positive basis]
A subset $P_1 \in \R^{n \times m}, 1 \leq m \leq n$ of a  basis $\B_n$ of $\R^n$ is called a \emph{sub-basis of a subspace $L_1$ in $\R^n$} if $\spann(P_1)=L_1.$ \\
A subset  $P_2 \in \R^{n \times r}, 2 \leq r \leq s,$ of a positive basis $\D_{n,s}$ is called a \emph{sub-positive basis of a subspace $L_2$ in $\R^n$} if $\pspan(P_2)=L_2.$ 
\end{definition}

A sub-positive basis of a subspace $L$ in $\R^n$  ($1 \leq \dim(L) \leq n$) will be denoted $\D^n_{m,r}$ where $m=\dim(L),$ $r$ is the  size of the sub-positive basis. When $m=n$ (that is the case when the subspace is $\R^n$ itself), we omit the superscript and simply write $\D_{n,s}$. When the sub-positive basis is minimal size, we may omit the size $r$ in the subscript, and simply write $\D^n_{m}.$

\begin{example}
Let $\D_{3,5}=\bbm 1&0&-\frac{1}{\sqrt{2}}&0&0\\0&1&-\frac{1}{\sqrt{2}}&0&0\\0&0&0&1&-1 \ebm.$ It can be proved that $\D_{3,5}$ is a positive basis of $\R^3.$ Also, we have that $$\D^3_{2}=\bbm 1&0&-\frac{1}{\sqrt{2}}\\0&1&-\frac{1}{\sqrt{2}}\\0&0&0 \ebm$$ is a sub-positive basis of $L_1=\{x \in \R^3:x_3=0\}$ and $$\mathring{\D}^3_1=\bbm0&0\\0&0\\1&-1\ebm$$ is a sub-positive basis of $L_2=\{x \in \R^3:x_1=x_2=0\}.$
\end{example}

\begin{definition}[Critical set, critical vectors and complete critical set] \label{def:criticalset}
Let $\D_{n,s}$ be a positive basis of $\R^n.$ Let $C$ be a subset  of $\R^n.$ We say $C$ is a \emph{critical set of $\D_{n,s}$} if
\begin{align} \label{eq:criticalsetC}
    \pspan \left ( (\D_{n,s} \setminus{\{d\}}) \cup C \right ) &\neq \R^n
\end{align}
 for each $d \in \D_{n,s}$. Elements of $C$ are called \emph{critical vectors}. The \emph{complete critical set} is denoted $C(\D_{n,s})$ and contains all critical set $C$ satisfying \eqref{eq:criticalsetC}.
\end{definition}

Note that  $\zero \in \R^n$ is a critical vector for every positive basis in $\R^n.$ The following example provides the critical set of a minimal positive basis in $\R^2.$ It is proved in \cite{Romanowicz1987} that for a minimal positive basis $\D_n$ in $\R^n, n \geq 2,$ we have
\begin{align}\label{eq:criticalsetminpb}
    C(\D_{n})&=-\bigcup_{i \neq j} \pspan \left ( \D_{n,s} \setminus {\{d_i,d_j\}} \right ).
\end{align}
Also, it is proved  in \cite{Romanowicz1987} that a maximal positive basis $\D_{n,2n}, n \geq 1,$ has the property 
\begin{align}\label{eq:cjmaximalpb}
    C(\D_{n,2n})&=\{ \zero \}. 
\end{align}

\begin{example}
Consider the minimal (optimal) positive basis 
\begin{align*}
    \mathring{\D}_2&=\bbm d_1&d_2&d_3 \ebm =\bbm 1 &-\frac{1}{2}& -\frac{1}{2} \\0 &\frac{\sqrt{3}}{2}&-\frac{\sqrt{3}}{2}\ebm.
\end{align*}
It follows from Equation \eqref{eq:criticalsetminpb} that $$C(\mathring{\D}_2)=-\pspan(d_1) \cup -\pspan(d_2) \cup -\pspan(d_3).$$
%Figure \ref{fig:criticalset} illustrates $C(\mathring{\D}_2).$
%\begin{figure}[ht]
%\centering
%\captionsetup{justification=centering}
%\includegraphics[width=0.8\textwidth]{CompleteCriticalSet.eps}
%\caption{The complete critical set $C(\mathring{\D}_2)$}
%\label{fig:criticalset}
%\end{figure}
\end{example}

The next theorem  provides the general structure of a non-minimal positive basis in $\R^n.$
\begin{theorem}[Structure of a non-minimal positive basis \textup{\cite[Theorem 1]{Romanowicz1987}}] \label{thm:romanothm1}
Let $n\geq 2$ and $s \geq n+2$. The set of $s$ vectors in $\R^n$, $D_{n,s}$,  is a positive basis of $\R^n$ if and only if $D_{n,s}$ admits the partition
\begin{align}\label{eq:structurerom}
    D_{n,s}&=\D^n_{m_1} \cup (\D^n_{m_2} \oplus c_1) \cup \dots \cup (\D^n_{m_q} \oplus c_{q-1})
\end{align}
 where $\D^n_{m_1}, \dots, \D^n_{m_q}$ are minimal sub-positive positive bases of subspaces $L_1, \dots, L_q$ of $\R^n,$ $\R^n=L_1 \oplus L_2 \oplus \dots \oplus L_q, L_i \cap L_j=\{\zero \}$ for $i \neq j,$ $1 \leq \dim L_i \leq n-1,$ and $c_j$ ($j \in \{1, \dots, q-1\}$) is a critical vector of the sub-positive basis $\D^n_{M_j, M_j+j}$ of the  subspace $L_1 \oplus \dots \oplus L_j$ of $\R^n$,  where
 $$\D^n_{M_1,M_1+1}=\D^n_{m_1}$$ and $$\D^n_{M_j}=\D^n_{m_1} \cup (\D^n_{m_2} \oplus c_1) \cup \dots \cup (\D^n_{m_j} \oplus c_{j-1})$$ for all $j \in \{2, \dots, q\}, q\geq 2.$
 \end{theorem}
 Let us provide an example to clarify the meaning of the previous theorem.
 \begin{example}
 Let us consider sets of 5 vectors in $\R^3.$ First, the set
 \begin{align*}
 D_{3,5}&=\bbm d_1&d_2&d_3&d_4&d_5 \ebm \\
 %&=\diag \left ( \mathring{\D}_{2}, \mathring{\D}_{1} \oplus \zero \right ) \\
 %&=\diag \left ( \mathring{\D}_{2}, \mathring{\D}_{1} \right ) \\
 &=\bbm 1 &-\frac{1}{2}& -\frac{1}{2}& 0&0\\0&\frac{\sqrt{3}}{2}&-\frac{\sqrt{3}}{2}&0&0\\0&0&0&1&-1  \ebm 
 \end{align*}
 is a positive basis of $\R^3,$ since $\mathring{\D}_2^3=\bbm d_1&d_2&d_3 \ebm$ is a minimal (optimal) sub-positive basis of the subspace $L_2=\{x \in \R^3:x_3=0\},$ the set $\mathring{\D}_1^3=\bbm d_4&d_5\ebm $ is a minimal (optimal) positive basis of the subspace $L_1=\{x \in \R^3:x_1=x_2=0\},$ $\R^3= L_2 \oplus L_1,  L_2 \cap L_1=\zero$ and $\zero$ is a critical vector for every positive basis. Hence, $D_{3,5}$ admits the partition $$D_{3,5}=\mathring{\D}_2^3 \cup \mathring{\D}_{1}^3.$$
 
 Second, the set
 \begin{align*}
     D'_{3,5}&=\bbm d_1&d_2&d_3&d_4'&d_5'\ebm  \\
     &=\bbm 1 &-\frac{1}{2}& -\frac{1}{2}& -1&-1\\0&\frac{\sqrt{3}}{2}&-\frac{\sqrt{3}}{2}&0&0\\0&0&0&1&-1  \ebm 
 \end{align*}
 is a positive basis of $\R^3.$ Notice $c=\bbm -1&0&0 \ebm^\top$ is a critical vector of the sub-positive basis $\mathring{\D}_2^3.$ Hence, $D'_{3,5}$ can be written as the partition
 \begin{align*}
     D'_{3,5}&=\mathring{\D}^3_{2} \cup  \left (\mathring{\D}^3_{1} \oplus c \right ).
 \end{align*}
 \end{example}
 
 The following lemma shows that the number of sub-positive bases in a partition of a positive basis given in \eqref{eq:structurerom} is $q=s-n.$ 
 \begin{lemma}
 Let $\D_{n,s}$ be a non-minimal  positive basis in $\R^n.$ Then the number of sub-positive bases $\D_{m_j}$ in the partition \eqref{eq:structurerom} is $$q=s-n.$$
 \end{lemma}
 \begin{proof}
 We have
\begin{align*}
    \sum_{j=1}^q (m_j+1)&=s,\\
    \sum_{j=1}^q m_j&=n,
\end{align*}
where $1 \leq m_j \leq n-1, n+1<s \leq 2n.$ Hence,
\begin{align*}
s-n&=\left ( \sum_{j=1}^q (m_j+1) \right )-n=q + \left (\sum_{j=1}^q m_j \right )-n = q.
\end{align*} \qed
 \end{proof}

 In the remainder of this paper, we focus on the positive bases that can be partitioned  such that  the  critical vectors $c_j$ are all equal to  $\zero.$ We will explain  why we decided to focus on this specific type of positive bases in the following pages.

\begin{definition}[The set $\Omega$]
 The positive basis  $\D_{n,s}$ of $\R^n$ is in the set $\Omega$ if  $\D_{n,s}$  admits the partition
 \begin{align}\label{eq:omega}
 \D_{n,s}&=\D^n_{m_1} \cup \D^n_{m_2} \cup \dots \cup \D^n_{m_{s-n}}
\end{align}
 where $\D^n_{m_1}, \dots, \D^n_{m_{s-n}}$ are minimal sub-positive positive bases of  the subspaces $L_1, \dots, L_{s-n}$  (respectively), $\R^n=L_1 \oplus L_2 \oplus \dots \oplus L_{s-n},$ $1 \leq \dim L_i=m_i \leq n,$ for all $i \in \{1, \dots,s-n\}$ and  such that $L_i \cap L_j=\{\zero \}$ for $i \neq j$ whenever $s-n \geq 2.$
 \end{definition}
 
 Note that the set $\Omega$ contains positive bases of all sizes: minimal, intermediate and maximal. Also, all maximal positive bases and minimal positive bases of $\R^n$ are contained in $\Omega.$ That is $\Omega=\Pe$ whenever $s=n+1$ or $s=2n.$ However, when $n+1<s<2n,$ we have that $\Omega$ is a proper non-empty subset of $\Pe.$
  The next corollary follows directly from Theorem \ref{thm:romanothm1}. It describes the structure of a basis $\B_n$ in $\R^n$  contained in a non-minimal   positive basis $\D_{n,s}$ contained in $\Omega$.
 \begin{corollary} \label{cor:Bnstructure}
 Let $\D_{n,s}$ be a positive basis  of $\R^n$ contained in $\Omega.$ Let $\B_n$ be any basis  of $\R^n$ contained in $\D_{n,s}.$ Then $\B_n$ admits the partition
 \begin{align*}
     \B_n&=\B^n_{m_1} \cup \B^n_{m_2} \cup \dots \cup \B^n_{m_{s-n}} 
 \end{align*}
 where $\B^n_{m_i} \in \R^{n \times m_i}$ is a sub-basis of the subspace $L_i$ in $\R^n$ for all $i \in \{1, \dots, s-n\}$, $\R^n=L_1 \oplus \dots \oplus L_{s-n},$ $1 \leq \dim L_i=m_i \leq n,$ $L_i \cap L_j=\{\zero\}$ for $i \neq j$ (whenever $s-n \geq 2$) and such that 
 \begin{align*}
     \B^n_{m_i} &\subset \D^n_{m_i}, \\
 \end{align*}
 for all $i \in \{1, \dots, s-n\}.$
 \end{corollary}
 
 In the previous corollary, note that a sub-positive basis $\D^n_{m_i}$ contains $m_i+1$  sub-bases of $L_i$ (this follows from Lemma \ref{lem:Removing}). Forming a set by picking one sub-basis from each subspace $L_i$ always forms a basis of $\R^n.$ Hence, the number of bases of $\R^n$ contained in a positive basis in $\Omega$ is 
 \begin{align} \label{eq:numberofbases}
 &\prod_{i=1}^{s-n}(m_i+1).
 \end{align}
 
 For example, we obtain that there are $n+1$ bases of $\R^n$ in a minimal positive basis and  $2^n$ bases of $\R^n$ in a maximal positive basis. This agrees with Proposition 14 and Proposition 15 in \cite{hare2020cm}.
 We conclude this section by showing that Algorithm \ref{alg:cmpbasis} developed in \cite{hare2020cm} can be simplified when $\D_{n,s}$ is a positive basis of $\R^n$ contained in $\Omega.$ 
 
 \begin{proposition}[Property of a positive basis of $\R^n$ contained in $\Omega$] \label{prop:propomega}
 Let $\D_{n,s}$ be a positive basis of $\R^n$ contained in $\Omega.$ Let $\B_n$ be a basis of $\R^n$ contained in $\D_{n,s}.$ Let $u_{\B_n}$ be the unit vectors such that $u_{\B_n}^\top \B_n=\gamma_{\B_n} \one^\top$ (where  $\gamma_{\B_n}$ is defined as in Algorithm \ref{algo:computecm} Step (1.1)). Then
 $$ u_{\B_n}^\top d \leq 0$$ for all vectors $d \in \D_{n,s} \setminus \{\B_n\}$. Consequently, Steps (1.3),(1.4) in Algorithm \ref{algo:computecm} can be omitted  and Steps (2.1),(2.2) become (respectively)
 \begin{align}
     \cm(\D_{n,s})&=\min_{\B_n \subset \D_{n,s}} \gamma_{\B_n},  \label{eq:newcm}\\
     \V(\D_{n,s})&=\{u_{\B_n}:\gamma_{\B_n}=\cm(\D_{n,s})\}. \label{eq:newV}
 \end{align}
 \end{proposition}
 \begin{proof}
 Let $\B_n$ be any basis of $\R^n$ contained in $\D_{n,s}.$ By Corollary \ref{cor:Bnstructure}, $\B_n$ can be written as $\B_n=\B^n_{m_1} \cup \dots \cup \B^n_{m_{s-n}}$  where $\B^n_{m_i}$ is contained in the sub-positive basis $\D^n_{m_i}$ of the subspace $L_i$ in $\R^n.$   Let $u_{\B_n}$ be the unit vector defined in  Step (1.2) of Algorithm \ref{alg:cmpbasis}.  Consider the minimal sub-positive basis $\D^n_{m_i}$  for some $i \in \{1, \dots, s-n\}.$ We know that there is only one vector in the set $\D^n_{m_i}\setminus \{\B^n_{m_i}\}.$ Denote this vector by $d$. We know that the projection of $u_{\B_n}$ onto $L_i,$ denoted $\proj_{L_i} u_{\B_n}$ has dot product $$ \left ( \proj_{L_i} u_{\B_n} \right )^\top d<0$$ whenever  $\proj_{L_i} u_{\B_n} \neq \zero$ by Lemma \ref{lem:dotprodneg}, and it is equal to zero when  $\proj_{L_i} u_{\B_n}=\zero.$  We obtain
 $$u_{\B_n}^\top d=\left ( \proj_{L_i} u_{\B_n} \right )^\top d \leq 0.$$
 Therefore, $u_{\B_n}^\top d \leq 0$ for all $d \in \D_{n,s} \setminus \{\B_n\}.$

 As mentioned in Section \ref{sec:prel}, we know that  $\gamma_{\B_n}>0$ for any basis $\B_n$ contained in $\D_{n,s}.$ It follows that $\gamma_{\B_n}$ is the maximum value of the dot product between $u_{\B_n}$ and the vectors in $\D_{n,s}.$  Therefore, we may delete Steps (1.3), (1.4) and simply set $\cm(\D_{n,s})$ and $\V(\D_{n,s})$ as in \eqref{eq:newcm} and \eqref{eq:newV} respectively. \qed
 \end{proof} 
 In the sequel, we consider positive bases in $\Omega.$  Let us begin by finding the structure of an optimal  non-minimal positive basis  in $\R^3$ over $\Omega.$
 
 %==================================================================================
 \subsection{Optimality in $\R^3$} \label{sec:optR3}
 In $\R^3,$ there is only one possible  intermediate size: 5. Let $\D_{3,5}$ be a positive basis contained in $\Omega.$  Then $\D_{3,5}$ admits the partition  
 \begin{align*}
     \D_{3,5}=\mathring{\D}^3_1 \cup \D^3_2 
 \end{align*}
 where $\mathring{\D}^3_1$ is a minimal  (optimal) sub-positive basis of  a subspace $L_1$ in $\R^3$ and $\D^3_2$ is a sub-positive basis of a subspace $L_2$ in $\R^3$ such that $L_1 \oplus L_2 =\R^3,$ $L_1 \cap L_2 = \{ \zero \}.$ Note that this is the only possible partition that includes at least two  sub-positive bases. Also, $\D_{3,5}$ cannot contain more than two sub-positive bases. We can realign the positive basis $\D_{3,5}$ so that $\mathring{\D}^3_1$ is  equal to
 \begin{align*}
     \mathring{\D}^3_1&=\bbm 0&0\\0&0\\1&-1 \ebm.
 \end{align*}
 This realignnment does not affect the geometry of the positive basis in the sense that it does not change the values of the Gram matrix $\G(\B_3).$ Hence, it does not affect the cosine measure.  It follows that any basis $\B_3$ of $\R^3$ contained in $\D_{3,5}$  admits the partition
 \begin{align*}
     \B_3=\B^3_1 \cup \B^3_2
 \end{align*}
 where $\B^3_1=\pm \bbm 0&0&1 \ebm^\top$ and $\B^3_2$ is a sub-basis of $L_2$ contained in $\D^3_2.$
 We next show that the two subspaces $L_1$ and $L_2$  of $\R^3$ must be orthogonal for $\D_{3,5}$ to be optimal over $\Omega.$ First, a lemma is introduced. It will be used in the in the main proposition of this sub-section.
 \begin{lemma} \label{lem:grandsumG3inverse}
Let $\D_{3,5}$ be a positive basis in $\R^3$ contained in $\Omega$. Let $\B_3$ be a basis of $\R^3$ contained in $\A(\D_{3,5})$ and admitting the partition $\B_3=\B_2^3 \cup \B_1^3.$ Let $v=(\B_3^2)^\top \B_1^3 \in \R^2.$ Then 
 \begin{align*}
     \one^\top \G(\B_3)^{-1} \one&=\one^\top \G(\B_2^3)^{-1} \one+c \left ( \one^\top \bbm -\G(\B_2^3)^{-1}v\\1\ebm\right )^2
 \end{align*}
 where $c=\left (1-v^\top \G(\B_2^3)^{-1}v \right )^{-1}.$ Moreover, 
 \begin{align}
     \one^\top \G(\B_3) \one &\geq \one^\top \G(\B_2^3)^{-1} \one+1
 \end{align}
 with equality if and only $v=\zero.$
 \end{lemma}
 \begin{proof}
  To make notation tighter, let $\G(\B_3)=\G_3$ and $\G(\B^3_2)=\G_2.$ By Lemma \ref{lem:inverseblock}, the inverse of $\G_3$ is
\begin{align*}
    \G_3^{-1}&=\bbm \G_2^{-1}+\G_2^{-1}v\left ( 1-v^\top \G_2^{-1}v\right )^{-1}v^\top \G_2^{-1}&\quad -\G_2^{-1}v\left (1-v^\top \G_2^{-1}v \right )^{-1}\\-\left (1-v^\top \G_2^{-1}v \right )^{-1} v^\top \G_2^{-1}&\left (1-v^\top \G_2^{-1}v \right )^{-1} \ebm.
\end{align*}
Once again, to make notation tighter, let $c=\left (1-v^\top \G_{2}^{-1}v\right )^{-1} \in \R.$  We obtain
\begin{align*}
    \G_3^{-1}&=\bbm \G_2^{-1}&\zero\\ \zero&0\ebm +c \bbm \G_2^{-1}vv^\top\G_2^{-1}&-\G_2^{-1}v\\(-\G_2^{-1}v)^\top&1 \ebm \\
    &=\bbm \G_2^{-1}&\zero \\ \zero&0\ebm+c \bbm-\G_2^{-1}v\\1 \ebm \bbm (-\G_2^{-1}v)^\top&1\ebm.
\end{align*}
It follows that the grand sum of $\G_3^{-1}$ is
\begin{align}\label{eq:grandsum}
    \one^\top \G_3^{-1} \one&=\one^\top \G_{2}^{-1}\one +c \left (\one^\top \bbm-\G_2^{-1}v\\1 \ebm \right )^2.
\end{align}
The second term in the previous equation is now investigated. We know that $\G_3^{-1}$ is a positive definite matrix by Lemma \ref{lem:inversepdispd}. Since $c$ is a principal submatrix of a positive definite matrix, $c$ is positive definite by Lemma \ref{lem:principalsubmatrix}. It follows that $c \geq 1$ with equality if and only if  $v=\zero.$ Now, we consider three cases.  If $-G_2^{-1}v=\zero,$ this means that $v_1+v_2=0.$ If $v_1=v_2=0,$ then 
$$c \left (\one^\top \bbm-\G_2^{-1}v\\1 \ebm \right )^2=1.$$
If $v_1+v_2=0$ and $v_1=-v_2\neq 0,$ then we get $c>1$ and so $$c \left ( \one^\top \bbm-\G_2^{-1}v\\1 \ebm\right )^2>1.$$ Lastly, we show that $-\G_2^{-1}v<0$ is not possible. Suppose $-\G_2^{-1}v<0$. Then $-\G_2^{-1}(-v)>0.$ Therefore, the basis form with $\B_2^3$ and $-\B_1^3,$ say $\widetilde{\B}_3$, has a grand sum $\one^\top \G(\widetilde{\B}_3)^{-1} \one $ strictly greater than $\one^\top \G(\B_3)^{-1} \one.$ This is a contradiction to the assumption that $\one^\top\G(\B_3)^{-1} \one$ is the maximum  for all possible bases of $\R^3$ contained in $\D_{3,5}$  since it is in $\mathbf{A}(\D_{3,5}).$ Therefore, we must have $-\G_2^{-1}v>0$ and it follows that $$c \left ( \one^\top \bbm-\G_2^{-1}v\\1 \ebm\right )^2>1.$$
Therefore, the second term in Equation \eqref{eq:grandsum} is minimal and equal to 1 if and only if $v=\zero.$ \qed
 \end{proof}

 \begin{proposition}[Orthogonality of the subspaces] \label{prop:orthoR3}
 Let $\D_{3,5}= \mathring{\D}^3_{1} \cup \D^3_{2}$ be a positive basis of size 5 in $\R^3$ where $\mathring{\D}^3_{1}$ is a sub-positive basis of the subspace $L_1$ and $\D^3_2$ is a sub-positive basis of the subspace $L_2,$ $L_1 \oplus L_2=\R^3$ such that $L_1 \cap L_2 =\{\zero\}.$  If $\D_{3,5}$ is optimal over $\Omega$, then $$L_1 \perp L_2.$$
 \end{proposition}
 \begin{proof}
 By contradiction. Suppose $\D_{3,5}$ is optimal and that the two subspaces are not orthogonal to each other. Let $\B_3$ be a basis of $\R^3$ in $\A(\D_{3,5}).$ We know that $\B_3$ admits a partition 
 \begin{align*}
     \B_3=\B^3_2 \cup \B^3_1,
 \end{align*}
 where $\B_2^3$ is contained in $\D_2^3$ and $\B_1^3$ is contained in $\mathring{\D}_1^3.$ By Lemma \ref{lem:grandsumG3inverse}, we conclude that the only possible way that $\D_{3,5}$ is optimal is to have $\one^\top \G(\B_2^3)^{-1} \one$ strictly less than the  best value for a positive basis where both subspaces are orthogonal since orthogonality of the subspaces decreases the grand sum $\one^\top \G(\B_3)^{-1} \one$ for a fix value of $\one^\top\G(\B_2^3)^{-1} \one$. We now show that it is not possible to obtain a value of $\one^\top \G(\B_2^3)^{-1} \one$ strictly less than  the value obtain when $\B^3_2$ is picked from an optimal sub-positive basis $\mathring{\D}^3_2,$ which is $2^2=4.$ 
By contradiction. Suppose that $\one^\top \G(\B_2^3)^{-1}\one <4.$ This means that there exists a sub-basis of $L_2$ contained in $\D^3_2$, say $\widetilde{\B}^3_2,$ such that $$\one^\top \G(\widetilde{\B}^3_2)^{-1} \one >4>\one^\top \G(\B_2^3)^{-1} \one.$$
Form $\widetilde{\B}_3$ by choosing $\widetilde{\B}^3_2$ and a vector $d$ contained in $\mathring{\D}^3_1$ such that $$c \left (\one^\top \bbm-\G(\widetilde{\B}_2^3)^{-1}v\\1 \ebm \right )^2>1,$$ where $v=(\widetilde{\B}_2^3)^\top d \in \R^2.$  Since $\B_3$ is in $\A(\D_{3,5}),$ it maximizes the grand sum $\one^\top \G(\cdot)^{-1} \one$ for all positive bases of $\R^3$ contained in $\D_{3,5}.$ Hence,
\begin{align*}
    \one^\top \G(\widetilde{\B}_3) \one &\leq \one^\top \G(\B_3)^{-1} \one.
\end{align*}
Let  $\D'_{3,5}=\diag (\mathring{\D}^3_2,\mathring{\D}^3_1)$ and $\B'_3 \in \A(\D'_{3,5})$. Since all terms in \eqref{eq:grandsum} for $\one^\top \G(\B'_3)^{-1} \one$ are strictly less than all corresponding terms in   \eqref{eq:grandsum} for $\one^\top \G(\widetilde{\B}_3)^{-1} \one,$ we obtain
$$\one^\top \G(\B'_3)^{-1} \one < \one^\top \G(\widetilde{\B}_3)^{-1} \one \leq \one^\top \G(\B_3)^{-1} \one.$$
By Proposition \ref{prop:propomega}, this means that $\cm(\D'_{3,5})>\cm(\D_{3,5}).$
This is a contradiction to the assumption that $\D_{3,5}$ is optimal.

Therefore, if $\D_{3,5}$ is optimal, then  the two  subspaces $L_1$ and $L_2$ must be orthogonal. \qed
 \end{proof}
 \begin{theorem}
 Let $\D_{3,5}=\diag  \left ( \mathring{\D}_2, \mathring{\D}_1 \right ).$ Then $\D_{3,5}$ is optimal  over $\Omega.$
 \end{theorem}
\begin{proof}
 Let $\D'_{3,5}=\D^3_2 \cup \mathring{\D}^3_1$ be an optimal positive basis over $\Omega$ where $\D_2^3$ is a sub-positive basis of the subspace $L_2$ and $\D^3_1$ is a sub-positive basis of the subspace $L_1$ in $\R^3$.  By Proposition \ref{prop:orthoR3}, $(\D^3_2)^\top \D^3_1=\zero \in \R^3.$ Let $\B_3=\B_2^3 \cup \B_1^3$ be a basis of $\R^3$ in $\A(\D'_{3,5}).$ Realigning the positive basis $\D'_{3,5}$ if necessary, the  Gram matrix of $\B_3$ is
\begin{align*}
    \G(\B_3)=\diag (\G(\B_2), 1),
\end{align*}
where $\B_2$ is a basis of $\R^2.$
The inverse of $\G(\B_3)$ is
\begin{align*}
    \G(\B_3)^{-1}&=\diag (\G(\B_2)^{-1}, 1).
\end{align*}
The cosine measure  is given by 
\begin{align*}
    \cm(\D'_{3,5})=\frac{1}{\sqrt{\one^\top \G(\B_2)^{-1} \one+ 1}}
\end{align*}
The minimal value of $\one^\top \G(\B_2)^{-1}\one$ is obtained if and only if $\B_2$ is picked from an optimal sub-positive  basis $\mathring{\D}^3_2$ of $\R^2.$ Hence, $\D_2^3=\mathring{\D}_2^3$ and we get that $\cm(\D'_{3,5})=\cm(\D_{3,5}).$ Therefore, $\D_{3,5}$ is optimal. \qed
\end{proof}
The following figure illustrates an optimal positive basis of $\R^3$ over $\Omega$  for each possible size ($s=4,5,6$). Note that $\mathring{\D}_3$ and $\mathring{\D}_{3,6}$ are also optimal over $\Pe$ since $\Omega=\Pe$ whenever $s \in \{n+1, 2n\}.$
\begin{figure}[H]
\centering
\includegraphics[width=0.85\textwidth]{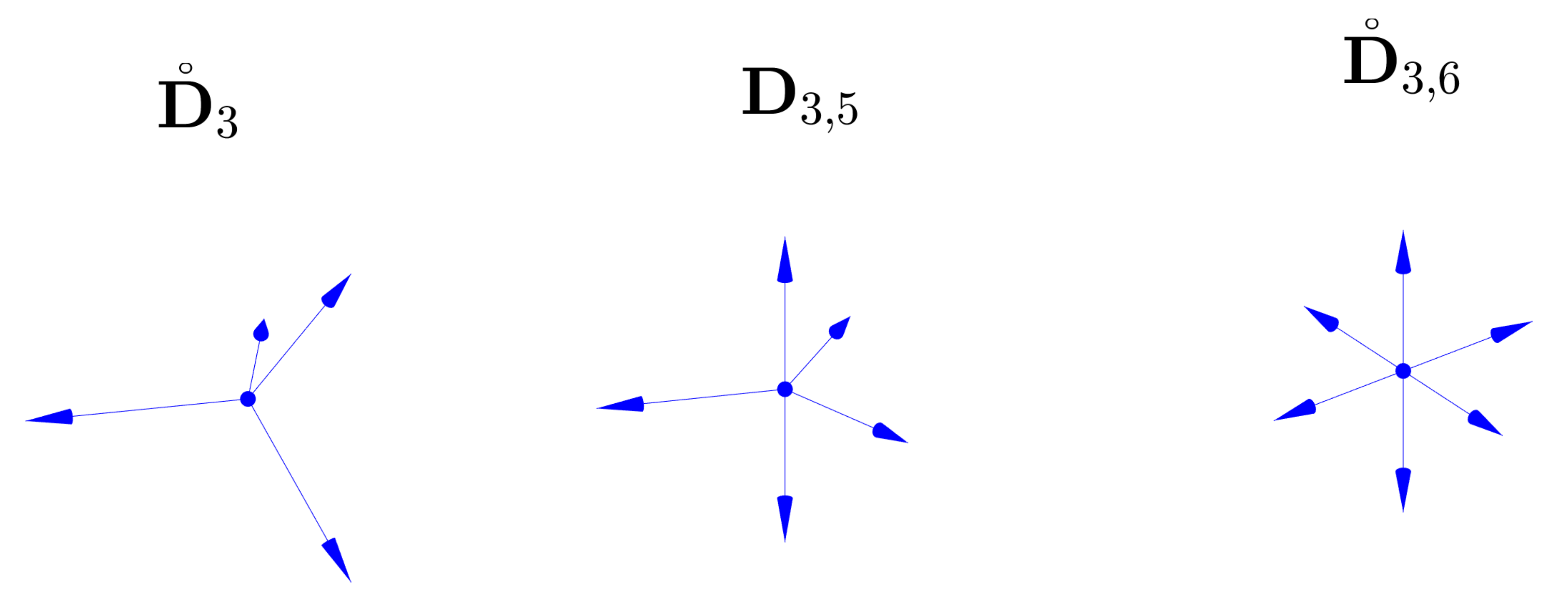}
\caption{An optimal positive basis of $\R^3$ over $\Omega$ for each possible size $s$.}
 \label{fig:Dstar3s}
\end{figure}
It is still unclear if  $\D_{3,5}=\diag (\mathring{\D}_2,\mathring{\D}_1 )$ is optimal over $\Pe.$ To prove that, it must be shown that there exists no positive basis of $\R^3$ with 5 vectors that admits a partition where not all the critical vectors are equal to zero that provides a greater cosine measure than $\cm(\D_{3,5}).$ A numerical experiment has been conducted and the results suggest that $\D_{3,5}$ is optimal over $\Pe,$ but no rigorous proof has been done yet. This topic is an obvious future research question to explore.

The next section investigates  optimality of a positive basis  in  a general space $\R^n.$
 \subsection{Optimality in $\R^n$} \label{sec:optRn}
 In the previous section, we have shown that the two subspaces of $\R^3$ must be orthogonal for a positive basis of size 5 to be optimal over $\Omega.$ The proof is relatively easy since one of the Gram matrix of a sub-basis is simply a scalar. This is not necessarily the case when considering a positive basis in higher dimensions. It is still open questions to determine if  all the subspaces must be pairwise orthogonal for a positive basis of $\R^n$ to be optimal over $\Omega$ and if all the critical vectors must be equal to zero for a positive basis of $\R^n$ to be optimal over $\Pe.$ Both these questions have been answered when $s=2n.$ Indeed, all the subspaces are pairwise orthogonal and all critical vectors are zero in a maximal positive basis. It seems reasonable to believe that it is also the case for intermediate positive bases. However, rigorously proving that it is the case (or proving that it is false) has still not be completed.   
 
 In the following section, we will restrict ourselves to  a proper subset of $\Omega$ that we named $\Omega^+$ and investigates the properties of an optimal positive basis over $\Omega^+.$ This is an important step to eventually  find a positive basis of intermediate size  over $\Pe.$ We will also see  that an optimal positive basis over $\Omega^+$ has a nice structure that makes it  easy to generate on a computer.

 \begin{definition}[The set $\Omega^+$]\label{def:omega+}
  Let $\D_{n,s}$ be a  positive basis in $\R^n.$ We say that $\D_{n,s}$ is in $\Omega^+$ if it is in $\Omega$ and all sub-positive bases $\D^n_{m_i}$   in a partition of $\D_{n,s}$ are pairwise orthogonal whenever $i \geq 2$. That is $(\D^n_{m_i})^\top \D^n_{m_j}=\zero$ for all $i \neq j, i,j \in \{1, \dots,, s-n\}$ whenever $s-n \geq 2.$
 \end{definition}
 
  Notice that, for positive bases of intermediate sizes, we have  $\Omega^+ \subset \Omega \subset \Pe.$ When $\D_{n,n+1}$ is a minimal positive basis, we have $\Omega^+=\Omega=\Pe.$ When $\D_{n,2n}$ is a maximal positive basis, we have $\Omega^+\subset \Omega=\Pe.$

 The next theorem presents two sufficient conditions for a positive basis in $\Omega^+$ to be optimal over $\Omega^+$. In  the next theorem, the notation $\rem(\frac{a}{b})$ is used to denote the remainder of the division $a/b,$ where $b$ is non-zero. Also, the notation $\vert I \vert,$ where $I$ is a finite index set, is used to represent the number of elements in $I.$
 \begin{theorem} \label{thm:optimalomega+}
 Let $\D_{n,s}$ be a positive basis of $\R^n$ in $\Omega^+.$ If the following two properties are satisfied, then $\D_{n,s}$ is optimal over $\Omega^+$.
 \begin{enumerate}[(i)]
 \item All the minimal sub-positive bases $\D^n_{m_i}$ involved in a partition of $\D_{n,s}$ are optimal, \label{item:2}
 \item the dimensions $m_i$ of the sub-positive bases $\D^n_{m_i}$ satisfy 
 \begin{align*}
    &m_j= \floor*{\frac{n}{s-n}}, \quad j \in J, &&m_k=\ceil*{\frac{n}{s-n}}, \quad k \in K,  
\end{align*}
where $J$ and $K$ are disjoint index set  such that 
\begin{align*}
&J \cup K=\{1, 2, \dots, s-n\}, &&\vert J \vert= s-n-\rem \left (\frac{n}{s-n} \right ), &&&\vert K \vert=\rem \left (\frac{n}{s-n} \right ).
\end{align*} \label{item:3}
 \end{enumerate}
 \end{theorem}
 \begin{proof}
 Suppose that $\D_{n,s}$ is a positive basis  of $\R^n$ in $\Omega^+$ such that properties  \textit{(i)}, \textit{(ii)} are satisfied and that $\D_{n,s}$ is not optimal. This means that there exists an optimal positive basis of size $s$ in $\R^n,$ say $\D'_{n,s},$ such that $\cm(\D'_{n,s})> \cm (\D_{n,s}).$   Let $\B_n \in \A(\D'_{n,s}).$ By Corollary \ref{cor:Bnstructure}, it follows that $\B_n$ admits the partition
 \begin{align*}
     \B_n=\B^n_{m_1} \cup \B^n_{m_{2}} \cup \dots \cup \B^n_{m_{s-n}}
 \end{align*}
 where $\B^n_{m_i} \in \R^{n \times m_i}$ is a sub-basis of the subspace $L_i$ of $\R^n$ for all $i \in \{1, \dots, s-n\}.$  Also, we know that  $\B^n_{m_i} \subset \D^n_{m_i}$ where $\D^n_{m_i}$ is a minimal sub-positive basis of $L_i$ for all $i \in \{1, \dots, s-n\}.$ Since $\D'_{n,s}$ is in $\Omega^+,$ all the subspaces $L_i$ are  orthogonal to each other. Hence, the Gram matrix of $\B_n$ is 
 \begin{align*}
     \G(\B_n)&=\diag( \G(\B^n_{m_1}), \dots, \G(\B^n_{m_{s-n}})).
 \end{align*}
 The inverse of $\G(\B_n)$ is
 \begin{align*}
 \G(\B_n)^{-1}&=\diag (\G(\B^n_{m_1})^{-1}, \dots, \G(\B^n_{m_{s-n}})^{-1}).
 \end{align*}
 Hence, the cosine measure of $\D'_{n,s}$ is given by 
 \begin{align*}
     \cm(\D'_{n,s})&=\frac{1}{\sqrt{\sum_{i=1}^{s-n}\one^\top \G(\B^n_{m_i})^{-1} \one }}.
 \end{align*}
 Since $\D'_{n,s}$ is optimal, we have that
 \begin{align*} 
 &\sum_{i=1}^{s-n} \one^\top \G(\B^n_{m_i})^{-1} \one
 \end{align*}
 is minimal. Hence, we must have that each  $\B^n_{m_i}$ is contained in an optimal minimal sub-positive basis $\mathring{\D}^n_{m_i}.$ So the sum in the previous equation is equal to
 \begin{align*}
     \sum_{i=1}^{s-n} \one^\top \G(\B^n_{m_i})^{-1} \one&=\sum_{i=1}^{s-n} m_i^2.
 \end{align*}
 The minimal possible value of $\sum_{i=1}^{s-n}m_i^2$ is obtained by solving the following optimization problem:
 \begin{align}\label{eq:optprobtosolve}
    \minimize & \sum_{i=1}^{s-n} m_i^2 \quad \st \sum_{i=1}^{s-n} m_i=n, \quad m_i \in \N.
\end{align}
The integer solution is given by
\begin{align*}
    &m_j= \floor*{\frac{n}{s-n}}, \quad j \in J, &&m_k=\ceil*{\frac{n}{s-n}}, \quad k \in K,  
\end{align*}
where $J$ and $K$ are disjoint index set such that 
\begin{align*}
&J \cup K=\{1, 2, \dots, s-n\}, &&\vert J \vert= s-n-\rem \left (\frac{n}{s-n} \right ), &&&\vert K \vert=\rem \left (\frac{n}{s-n} \right ).
\end{align*}
 But then we obtain that  
 \begin{align*}
     \cm(\D_{n,s})&=\cm(\D'_{n,s}).
 \end{align*}
 A contradiction. Therefore, a positive  basis of $\R^n$  satisfying Properties  \textit{(i)}, \textit{(ii)} must be optimal over $\Omega^+.$ \qed
 \end{proof}
 \begin{corollary}[The cosine measure of an optimal positive basis over $\Omega^+$]
 Let $\D_{n,s}$ be an optimal positive basis over $\Omega^+.$ Let $\mathtt{r}=\rem \left ( \frac{n}{s-n}\right ).$ Then
 \begin{align*}
     \cm(\D_{n,s})&=\frac{1}{\sqrt{(s-n-\mathtt{r})\floor*{\frac{n}{s-n}}^2+\mathtt{r}\ceil*{\frac{n}{s-n}}^2}}.
 \end{align*}
 \end{corollary}
In particular, when $\mathring{\D}_{n}$ is an optimal minimal positive basis over $\Omega^+(=\Pe)$, we obtain $$\cm({\mathring{\D}_n})=\frac{1}{n}.$$ This agrees with the value of the cosine measure provided in \cite[Theorem 1]{Naevdal2018} for a minimal positive basis to be optimal over $\Pe$. When $\D_{n,2n}$ is an optimal positive basis over $\Omega^+,$ we obtain $$\cm(\D_{n,2n})= \frac{1}{\sqrt{n}}.$$ Once again, this value agrees with the value provided in \cite[Theorem 2]{Naevdal2018} for a maximal positive basis to be optimal over $\Pe.$
Hence, for both minimal and maximal positive bases to be optimal over $\Pe,$ it is necessary that they are contained in $\Omega^+.$ This provides an argument to believe that this is also the case for intermediate positive bases. However, two facts remain to be proved (or disproved) in $\R^n$ before concluding that it is the case (or not). 

The following table provides the diagonal blocks contained in  an optimal (over $\Omega^+$) positive basis of $\R^n$ of the form $\D_{n,s}=\diag(\mathring{\D}_{m_1}, \dots, \mathring{\D}_{m_{s-n}}).$  The notation $(\mathring{\D}_{m_i})^k$ where $k$ is a positive integer means that the diagonal block $\mathring{\D}_{m_i}$ appears $k$ times as a diagonal entry in $\D_{n,s}.$
\renewcommand{\arraystretch}{1.2}

\begin{table}[H]
\centering
\scalebox{1}{
\begin{tabular}{|l|l|l|l|l|l|l|l|}\hline
$s \setminus \,n$&2&3&4&5&6&7&8\\
\hline
&&&&&&&\\[-1em]
3&$\mathring{\D}_{2}$&-&-&-&-&-&-\\\hline
&&&&&&&\\[-1em]
4&$(\mathring{\D}_{1})^2$&$\mathring{\D}_{3}$&-&-&-&-&-\\\hline
&&&&&&&\\[-1em]
5&-&$\mathring{\D}_{2},\mathring{\D}_{1}$&$\mathring{\D}_{4}$&-&-&-&-\\\hline
&&&&&&&\\[-1em]
6&-&$(\mathring{\D}_{1})^3$&$(\mathring{\D}_2)^2$&$\mathring{\D}_{5}$&-&-&-\\ \hline
&&&&&&&\\[-1em]
7&-&-&$\mathring{\D}_2,(\mathring{\D}_{1})^2$&$\mathring{\D}_2,\mathring{\D}_{3}$&$\mathring{\D}_{6}$&-&- \\
\hline
&&&&&&&\\[-1em]
8&-&-&$(\mathring{\D}_{1})^4$&$(\mathring{\D}_2)^2, \mathring{\D}_{1}$&$(\mathring{\D}_{3})^2$&$\mathring{\D}_{7}$&- \\
\hline
&&&&&&&\\[-1em]
9&-&-&-&$\mathring{\D}_2,(\mathring{\D}_{1})^3$&$(\mathring{\D}_2)^3$&$\mathring{\D}_{4},\mathring{\D}_{3}$&$\mathring{\D}_{8}$\\
\hline
&&&&&&&\\[-1em]
10&-&-&-&$(\mathring{\D}_{1})^5$&$(\mathring{\D}_2)^2,(\mathring{\D}_{1})^2$&$\mathring{\D}_{3},(\mathring{\D}_2)^2$&$(\mathring{\D}_{4})^2$ \\
\hline
&&&&&&&\\[-1em]
11&-&-&-&-&$\mathring{\D}_2,(\mathring{\D}_{1})^4$&$(\mathring{\D}_2)^3,\mathring{\D}_{1}$&$(\mathring{\D}_{3})^2,\mathring{\D}_2$ \\
\hline
&&&&&&&\\[-1em]
12&-&-&-&-&$(\mathring{\D}_{1})^6$&$(\mathring{\D}_2)^2,(\mathring{\D}_{1})^3$&$(\mathring{\D}_2)^4$ \\
\hline
&&&&&&&\\[-1em]
13&-&-&-&-&-&$\mathring{\D}_2,(\mathring{\D}_{1})^5$&$(\mathring{\D}_2)^3,(\mathring{\D}_{1})^2$ \\
\hline
&&&&&&&\\[-1em]
14&-&-&-&-&-&$(\mathring{\D}_{1})^7$&$(\mathring{\D}_2)^2,(\mathring{\D}_{1})^4$ \\
\hline
&&&&&&&\\[-1em]
15&-&-&-&-&-&-&$\mathring{\D}_2,(\mathring{\D}_{1})^6$ \\ 
\hline
&&&&&&&\\[-1em]
16&-&-&-&-&-&-&$(\mathring{\D}_{1})^8$ \\
\hline
%Chain Rule&\\\hline
\end{tabular}}  \caption{The diagonal blocks in  an optimal (over $\Omega^+$) positive basis $\D_{n,s}=\diag(\mathring{\D}_{m_1}, \dots, \mathring{\D}_{m_{s-n}})$} \label{table:dstarcm}
\end{table} 
A Matlab code is available on request to generate an optimal positive basis over $\Omega^+$ of any dimension $n$ and size $s.$ Note that each minimal sub-positive basis may be realigned in their respective subspace to include a specific vector of the subspace. The whole positive basis may also be realigned to include  a specific vector of $\R^n.$ These realignments do not affect the value of the cosine measure as it can be done by multiplying $\D_{n,s}$ with an orthonormal matrix.   A method to  accomplish these realignments  is provided in \cite{Jarry2019}.

\section{Conclusion} \label{sec:conc}
We have investigated the structure of intermediate positive bases in $\R^n$.  Using results from Romanowicz, we demonstrated that any positive basis, $\D_{n,s}$,  of $\R^n$ can be partitioned  into $s-n$ minimal sub-positive bases plus a critical vector. In this paper, we have focused on $\Omega$, the set of positive bases that can be written as a partition in which all critical vectors are zero.  We have shown that the algorithm developed in \cite{hare2020cm} to compute the cosine measure can be simplified for any positive basis contained in $\Omega$ (Proposition \ref{prop:propomega}).  When $\D_{n,s}$ is a positive basis contained in $\Omega,$ the number of bases of $\R^n$ contained in $\D_{n,s}$ has been identified \eqref{eq:numberofbases}. We have provided the structure of an optimal positive basis of intermediate size over $\Omega$ in $\R^3$, $\D_{3,5}$, and proved that the two subspaces must be orthogonal to each other.   We conjecture that this result also hold in $\R^n$.  In reviewing this result, note that the key requirement would be to extend Lemma \ref{lem:inverseblock} to include a broader class of positive bases.

Following this conjecture, we focused on $\Omega^+$, the set of positive bases that  an be written as a partition in which all critical vectors are zero and all sub-bases are orthogonal.  We have determine a characterization of the structure of an optimal positive basis of this type (Theorem \ref{thm:optimalomega+}) and therefore determined the optimal cosine measure for positive bases in $\Omega^+$.  It turns out that it is very simple and efficient to generate such a positive basis with software such as Matlab.  A Matlab code is available upon request.

In order to characterize optimality over the whole universe $\Pe$ for positive bases of intermediate size completely, two questions need to be answered.  First, must all critical vectors in a partition of an optimal positive basis be zero? Second, must all subspaces involved in a partition of an optimal positive basis be orthogonal?
We conjecture that the answer to both questions is yes and will further examine this in future research. 

%===============================================================================
% Model of figure 
%================================================================================
\begin{comment}

\begin{figure}[ht]
\centering
\includegraphics[height=6cm]{vonY.png}
\caption{$\von(\Y_\epsilon)$ as a function of $\epsilon$ for different $n$}
\label{Fig:VonY}
\end{figure}

\end{comment}

\normalsize
\bibliographystyle{siam}
\bibliography{Mybib}

\begin{thebibliography}{10}

\bibitem{Abramson2009}
{\sc M.~Abramson, C.~Audet, J.~Dennis, and S.~Le~Digabel}, {\em Orthomads: A
  deterministic mads instance with orthogonal directions}, SIAM Journal on
  Optimization, 20 (2009), pp.~948--966.

\bibitem{Audet2011}
{\sc C.~Audet}, {\em A short proof on the cardinality of maximal positive
  bases}, Optimization Letters, 5 (2011), pp.~191--194.

\bibitem{Audet2006}
{\sc C.~Audet and J.~Dennis}, {\em Mesh adaptive direct search algorithms for
  constrained optimization}, SIAM Journal on optimization, 17 (2006),
  pp.~188--217.

\bibitem{Audet2017}
{\sc C.~Audet and W.~Hare}, {\em Derivative-free and blackbox optimization},
  Springer Series in Operations Research and Financial Engineering, Springer,
  Cham, 2017.

\bibitem{Audet2014}
{\sc C.~Audet, A.~Ianni, S.~Le~Digabel, and C.~Tribes}, {\em Reducing the
  number of function evaluations in mesh adaptive direct search algorithms},
  SIAM Journal on Optimization, 24 (2014), pp.~621--642.

\bibitem{beck2020}
{\sc A.~Beck and N.~Hallak}, {\em On the convergence to stationary points of
  deterministic and randomized feasible descent directions methods}, SIAM
  Journal on Optimization, 30 (2020), pp.~56--79.

\bibitem{Conn2009}
{\sc A.~Conn, K.~Scheinberg, and L.~Vicente}, {\em Introduction to
  derivative-free optimization}, vol.~8, Siam, 2009.

\bibitem{Coope2001}
{\sc I.~Coope and C.~Price}, {\em On the convergence of grid-based methods for
  unconstrained optimization}, SIAM Journal on Optimization, 11 (2001),
  pp.~859--869.

\bibitem{Coope2002}
\leavevmode\vrule height 2pt depth -1.6pt width 23pt, {\em Positive bases in
  numerical optimization}, Computational Optimization and Applications, 21
  (2002), pp.~169--175.

\bibitem{Custodio2008}
{\sc A.~Cust{\'o}dio, J.~Dennis, and L.~Vicente}, {\em Using simplex gradients
  of nonsmooth functions in direct search methods}, IMA Journal of Numerical
  Analysis, 28 (2008), pp.~770--784.

\bibitem{Davis1954}
{\sc C.~Davis}, {\em Theory of positive linear dependence}, American Journal of
  Mathematics, 76 (1954), pp.~733--746.

\bibitem{hare2020cm}
{\sc W.~Hare and G.~Jarry-Bolduc}, {\em A deterministic algorithm to compute
  the cosine measure of a finite positive spanning set}, Optimization Letters,
  14 (2020), pp.~1305--1316.

\bibitem{Horn1990}
{\sc R.~Horn and C.~Johnson}, {\em Matrix analysis}, Cambridge university
  press, 1990.

\bibitem{Jarry2019}
{\sc G.~Jarry-Bolduc, P.~Nadeau, and S.~Singh}, {\em Uniform simplex of an
  arbitrary orientation}, Optimization Letters,  (2019), pp.~1--11.

\bibitem{Kelley2011}
{\sc C.~Kelley}, {\em Implicit filtering}, vol.~23, SIAM, 2011.

\bibitem{Kolda2003}
{\sc T.~Kolda, R.~Lewis, and V.~Torczon}, {\em Optimization by direct search:
  New perspectives on some classical and modern methods}, SIAM review, 45
  (2003), pp.~385--482.

\bibitem{Lewis1996}
{\sc R.~Lewis and V.~Torczon}, {\em Rank ordering and positive bases in pattern
  search algorithms.}, tech. rep., Institute for Computer Applications in
  Science and Engineering, Hampton VA, 1996.

\bibitem{Lu2002}
{\sc T.~Lu and S.~Shiou}, {\em Inverses of 2$\times$ 2 block matrices},
  Computers \& Mathematics with Applications, 43 (2002), pp.~119--129.

\bibitem{Mckinney1962}
{\sc R.~McKinney}, {\em Positive bases for linear spaces}, Transactions of the
  American Mathematical Society, 103 (1962), pp.~131--148.

\bibitem{Naevdal2018}
{\sc G.~N{\ae}vdal}, {\em Positive bases with maximal cosine measure},
  Optimization Letters,  (2018), pp.~1--8.

\bibitem{Regis2016}
{\sc R.~Regis}, {\em On the properties of positive spanning sets and positive
  bases}, Optimization and Engineering, 17 (2016), pp.~229--262.

\bibitem{regis2021}
{\sc R.~Regis}, {\em On the properties of the cosine measure and the uniform
  angle subspace}, Computational Optimization and Applications, 78 (2021),
  pp.~915--952.

\bibitem{Romanowicz1987}
{\sc Z.~Romanowicz}, {\em Geometric structure of positive bases in linear
  spaces}, Applicationes Mathematicae, 19 (1987), pp.~557--567.

\bibitem{Shephard1971}
{\sc G.~Shephard}, {\em Diagrams for positive bases}, Journal of the London
  Mathematical Society, 2 (1971), pp.~165--175.

\bibitem{Torczon1997}
{\sc V.~Torczon}, {\em On the convergence of pattern search algorithms}, SIAM
  Journal on optimization, 7 (1997), pp.~1--25.

\bibitem{van2013using}
{\sc B.~Van~Dyke and T.~Asaki}, {\em Using qr decomposition to obtain a new
  instance of mesh adaptive direct search with uniformly distributed polling
  directions}, Journal of Optimization Theory and Applications, 159 (2013),
  pp.~805--821.

\bibitem{Vaz2009}
{\sc A.~Vaz and L.~Vicente}, {\em Pswarm: a hybrid solver for linearly
  constrained global derivative-free optimization}, Optimization Methods \&
  Software, 24 (2009), pp.~669--685.

\end{thebibliography}
%\nocite{*}

%\appendix
%begin{appendices}
%\clearpage
%\pagenumbering{roman}
%\section*{Appendix}
%\lstinputlisting{cosmeasureps.m}
%\clearpage
%\lstinputlisting{pSpanCheck.m}
%\end{appendices}
\end{document}